\documentclass[12pt,draftcls,onecolumn]{IEEEtran}

\usepackage{amsfonts}

\IEEEoverridecommandlockouts                              
\usepackage{amsmath}
\usepackage{amssymb}

\title{On Sontag's Formula for the Input-to-State Practical Stabilization of Retarded Control-Affine Systems  }

\author{
P. Pepe {\it \ \ }
\thanks{
This work was supported in part by the Center of Excellence for
Research DEWS, the Italian MIUR Project PRIN 2009 and Atheneum Project RIA.
}
\thanks{P. Pepe is with
the Department of Information Engineering, Computer Science, and Mathematics, 
University of L'Aquila, L'Aquila, Italy,
{\tt\small pierdomenico.pepe@univaq.it}.}%
\thanks{
}%
}

\newtheorem{definition}{Definition}
\newtheorem{theorem}{Theorem}

\newtheorem{remark}{Remark}

\newtheorem{hypothesis}{Hypothesis}

\def\bn{\begin{eqnarray*}}
\def\be{\begin{eqnarray}}
\def\en{\end{eqnarray*}}
\def\ee{\end{eqnarray}}
\def\ba{\begin{array}}
\def\ea{\end{array}}
\def\bm{\left[\begin{array}}
\def\em{\end{array}\right]}

\begin{document}
\addtolength{\baselineskip}{-.08mm}

\maketitle
\thispagestyle{empty}   
\pagestyle{empty}       

\begin{abstract}
In this paper input-to-state practically stabilizing control laws for retarded,
control-affine, nonlinear systems with actuator disturbance are investigated.
The developed methodology is based on the Arstein's theory of control Liapunov
functions and related Sontag's formula, extended to retarded systems.
If the actuator disturbance is bounded, then the controller yields the solution
of the closed-loop system to achieve an arbitrarily fixed neighborhood of the
origin, by increasing a control tuning parameter. The considered systems can
present an arbitrary number of discrete as well as distributed time-delays, of
any size, as long as they are constant and, in general, known.
\noindent
{\it @ The extended version of this paper is under review in Systems \& Control Letters, Elsevier.}

\end{abstract}

\IEEEpeerreviewmaketitle

\section{Introduction}

Stabilization of nonlinear retarded systems is a topic which has
attracted many researchers in the last ten years. Many
contributions concerning the state feedback stabilization and the
input-output state feedback linearization of nonlinear time-delay
systems can be found in the literature (see, for instance,
\cite{kybGMP}\cite{GerManPeIJRNC}\cite{HuaGuanShi}\cite{Jankovic}\cite{Lien}\cite{MarquezMoog}\cite{OguchiIJC}\cite{PepeIJACSP}\cite{ZhangCheng}).
 The technique of control
Liapunov functions has been exploited to practically stabilize or
stabilize a large class of time-invariant time-delay systems in
affine form in \cite{Jankovic}, using Liapunov Razumikhin
functions. It is pointed out in \cite{Jankovic} that the Sontag's
formula cannot be directly applied to retarded systems by the use
of control Liapunov-Razumikhin functions, since the resulting
control input may well be discontinuous and unbounded when the
Razumikhin condition is not satisfied in a piece of the solution
trajectory. For this reason, in \cite{Jankovic} the domination
redesign control methodology is employed. It is also well known,
as far as the stability of time-delay systems is concerned, that
the Razumikhin method can be considered as a particular case of
the method of Liapunov–Krasovskii functionals (see
\cite{NEWKolmanovskii}, Section 4.8, p. 254). To our best
knowledge, first results concerning the use of Control
Liapunov-Krasovskii functionals (instead of Control Liapunov
Razumikhin functions) for building stabilizing control laws for
time-delay systems can be found in
\cite{MahboobiEsfanjaniNikravesh} and in \cite{KarJiangESAIM}. In
\cite{MahboobiEsfanjaniNikravesh} the Authors propose a predictive
control scheme with guaranteed closed-loop stability for
non-linear time-delay systems utilising a fixed class of control
Liapunov–Krasovskii functionals. In \cite{KarJiangESAIM} the
Authors prove the equivalence of the existence of a Control
Liapunov-Krasovskii functional and the stabilizability property,
for fully nonlinear, retarded systems. Moreover, stabilizability
is intended robust with respect to suitable disturbances,
un-effective when the state is zero. The problem of the robust
stabilization for a very large class of time-delay systems in
triangular form is extensively studied and solved. The reader is
referred to the recent monograph \cite{KarafyllisJiangBook} for
the topic of robust stability and stabilization of nonlinear
systems, in both the finite dimensional and the retarded case, in
continuous and discrete time.

In this paper we show how invariantly differentiable functionals
(see \cite{Kim}, \cite{KimBook}) can play an important role for
the input-to-state practical stabilization of retarded systems.
Invariantly differentiable functionals are very often used in the
literature concerning time-delay systems stability, though its
definition is not, according to us, very popular. The definition
of invariant differentiability was helpful for finding a
robustifying control law, yielding input-to-state stability with
respect to actuator disturbance, in \cite{PepeTAC09}. Here we show
that the same definition is very helpful in order to apply the
Sontag's formula to retarded, control-affine, nonlinear systems.
In details, the definition of invariant differentiability allows
to split the functional derivative (in the Driver's form, see
\cite{PepeJiangSCL} and references therein) into different parts.
Each part is suitably used in the Sontag's formula (i.e., in the
formulation of the terms $a$, $b$ in \cite{SontagUniversal}). In
this paper the Sontag's formula extended to retarded systems is
modified in the critical subsets of the infinite dimensional state
space where the Lipschitz property of the related feedback control
law may be lost. By this modification, the problem of non
Lipschitz feedback control law, which may cause non uniqueness or
even non existence of the solution, as well as implementation
difficulties, is solved. Then, a Sontag's input-to-state
stabilizing term (see \cite{SontagPrincipe}, \cite{PepeTAC09}) is
added to the control law, thus achieving the twofold result of
attenuation of the actuator disturbance and attenuation of the
bounded error due to the above modification of the Sontag's
formula. If the actuator disturbance is bounded, then an
arbitrarily fixed neighborhood of the origin is asymptotically
reached by the closed-loop solution, by increasing a tuning
parameter of the control law. The class of control affine systems
considered in this paper is more general than the one considered
in \cite{Jankovic}, since no hypothesis is introduced for the form
of the maps describing the dynamics (for instance, as far as the
presence of discrete or distributed delay terms is concerned).

 {\bf Notations}

\noindent $R$ denotes the set of real numbers, $R^{\star}$ denotes
the extended real line $[-\infty,+\infty]$, $R^+$ denotes the set
of non negative reals $[0,+\infty)$. The symbol $\vert \cdot
\vert$ stands for the Euclidean norm of a real vector, or the
induced Euclidean norm of a matrix. The essential supremum norm of
an essentially bounded function is indicated with the symbol
$\Vert \cdot \Vert_{\infty}$. A function $v:R^+\to R^m$, $m$
positive integer, is said to be {\it essentially bounded} if
$ess\sup_{t\ge 0}\vert v(t) \vert<+\infty$. For given times $0\le
T_1<T_2$, we indicate with $v_{[T_1,T_2)}:R^+\to R^m$ the function
given by $ v_{[T_1,T_2)}(t)= v(t) $ for all $t \in [T_1,T_2)$ and
$=0$ elsewhere. An input $v$ is said to be {\it locally
essentially bounded} if, for any $T>0$, $v_{[0,T)}$ is essentially
bounded. For a positive integer $n$, for a positive real $\Delta$
(maximum involved time-delay), ${\mathcal C}$ and ${\mathcal Q}$
denote the space of the continuous functions mapping $[-\Delta,0]$
into $R^n$ and the space of the bounded, continuous except at a
finite number of points, and right-continuous functions mapping
$[-\Delta,0)$ into $R^n$, respectively. For $\phi\in {\mathcal
C}$, $\phi_{[-\Delta,0)}$ is the function in ${\mathcal Q}$
defined as $\phi_{[-\Delta,0)}(\tau)=\phi(\tau)$, $\tau \in
[-\Delta,0)$. For a function $x:[-\Delta, c)\to R^n$, with $0<c\le
+\infty$, for any real $t\in [0,c)$, $x_t$ is the function in
${\mathcal C}$ defined as $x_{t}(\tau)=x(t+\tau)$, $\tau \in
[-\Delta,0]$. For a positive real $\delta$, $\phi\in {\mathcal
C}$, $I_{\delta}(\phi)=\{\psi \in {\mathcal C}: \Vert
\psi-\phi\Vert_{\infty}<\delta\}$, $\overline I_{\delta}(\phi)$ is
the closure of $I_{\delta}(\phi)$. For given positive integers
$n,m$, a map $f: \mathcal C\to R^{n\times m}$ is said to be:
completely continuous if it is continuous and takes closed bounded
subsets of ${\mathcal C}$ into bounded subsets of $R^{n\times m}$;
locally Lipschitz in ${\mathcal C}$ if, for any $\phi \in
{\mathcal C}$, there exist positive reals $\delta, \eta$ such
that, for any $\phi_1, \phi_2 \in \overline I_{\delta}(\phi)$, the
inequality $\vert f(\phi_1)-f(\phi_2)\vert \le \eta \Vert
\phi_1-\phi_2\Vert_{\infty}$ holds.

 Let us
here recall that a function $\gamma:R^+\to R^+$ is: of class
${\cal P}$ if it is continuous, zero at zero, and positive at any
positive real; of class ${\cal K}$ if it is of class ${\cal P}$
and strictly increasing; of class ${\cal K}_{\infty}$ if it is of
class ${\cal K}$ and it is unbounded; of class ${\cal L}$ if it is
continuous and it monotonically decreases to zero as its argument
tends to $+\infty$. A function $\beta:R^+\times R^+\to R^+$ is of
class ${\cal KL}$ if $\beta(\cdot, t)$ is of class ${\cal K}$ for
each $t\ge 0$ and $\beta(s,\cdot)$ is of class ${\cal L}$ for each
$s\ge 0$.  With $M_a$ (see \cite{PepeTAC09}) is indicated any
functional mapping ${\mathcal C}$ into $R^+$ such that, for some
${\mathcal K}_{\infty}$ functions $\gamma_{a}, \overline
\gamma_{a}$, the following inequalities hold
\begin{equation}\label{generaljiangcondition}
\gamma_{a}(\vert \phi(0) \vert) \le M_a(\phi) \le \overline
\gamma_{a}(\Vert \phi\Vert_{\infty}), \ \forall \ \phi \in
{\mathcal C}
\end{equation}
For example, the $\Vert \cdot \Vert_{M_2}$ norm, given by $\Vert
\phi \Vert_{M_2}= \left ( \vert
\phi(0)\vert^2+\int_{-\Delta}^0\vert \phi(\tau)\vert^2d\tau\right
)^{1\over 2}$, fulfills the conditions
(\ref{generaljiangcondition}) and thus is an $M_a$ functional.
RFDE stands for retarded functional differential equation, ISS
stands for input-to-state stability or input-to-state stable, ISpS
stands for input-to-state practical stability or input-to-state
practically stable, GAS stands for global asymptotic stability or
globally asymptotically stable. A system with an equilibrium at
zero is said $0$-GAS if the zero solution is GAS.

\section {Preliminaries}

Let us consider the following RFDE

\begin{eqnarray}\label{equazionediscontinua}&&
\dot x(t)=f(x_t)+g(x_t)u(t), \qquad t\ge 0,  \nonumber
\\ && x(\tau)=\xi_0(\tau),\qquad \tau \in [-\Delta,0], \qquad   \xi_0\in {\mathcal C},
\end{eqnarray}
where: $x(t)\in R^n$, $n$ is a positive integer; $\Delta>0$ is the
maximum involved time-delay; the maps $f:{\mathcal C}\to R^n$ and
$g:{\mathcal C}\to R^{n\times m}$ are completely continuous and
locally Lipschitz in ${\mathcal C}$; $m$ is a positive integer;
$u(t)\in R^m$ is the input signal, Lebesgue measurable and locally
essentially bounded.

The following definition of invariant differentiable functionals
is taken from \cite{KimBook}, see Definitions 2.2.1, 2.5.2 in
Chapter 2. The formalism used in \cite{KimBook} is here slightly
modified for the purpose of formalism uniformity over the paper.
For any given $x\in R^n$, $\phi \in {\mathcal Q}$ and for any
given continuous function ${\cal Y}:[0,\Delta]\to R^n$ with ${\cal
Y}(0)=x$, let $\psi^{(x,\phi,{\cal Y})}_h \in {\mathcal Q}$, $h\in
[0,\Delta)$, be defined as
\begin{eqnarray}\label{psixphiy}&&
\psi^{(x,\phi,{\cal Y})}_0=\phi; \nonumber \\ && for \ h>0, \qquad
\psi^{(x,\phi,{\cal Y})}_h (s) = \left \{
\begin{array}{c} \phi(s+h), \qquad \qquad s\in [-\Delta,-h);   \\
{\mathcal Y}(s+h), \qquad \qquad s\in [-h,0)\end{array}\right .
\end{eqnarray}
Let, for $\phi\in {\mathcal C}$, $h\in [0,\Delta)$, $\phi^h\in
{\mathcal C}$ be defined as follows
\begin{equation}
\phi^h(s)=\left \{\begin{array}{cc}\phi(s+h), & s\in
[-\Delta,-h)\\ \phi(0), & s\in [-h,0]\end{array}\right .
\end{equation}
\begin{definition}\label{invariantV}(see \cite{KimBook})
A functional $V:R^n\times{\mathcal Q}\to R^+$ is said to be
invariantly differentiable if, at any point $(x,\phi)\in R^n\times
{\mathcal Q}$: i) for any continuous function ${\cal
Y}:[0,\Delta]\to R^n$ with ${\cal Y}(0)=x$, there exists finite
the right-hand derivative $\left .\frac {\partial V\left ( x,
\psi^{(x,\phi,{\cal Y})}_h\right )}{\partial h}\right |_{h=0}$ and
such derivative is invariant with respect to the function ${\cal
Y}$;  ii) there exists finite the derivative $ \frac{\partial
V(x,\phi)}{\partial x}$;  iii) for any $z\in R^n$, for any
continuous function ${\cal Y}:[0,\Delta]\to R^n$ with ${\cal
Y}(0)=x$, for any $h\in [0,\Delta)$,
\begin{eqnarray}\label{differenziale}&& V\left (x+z,\psi^{(x,\phi,{\cal Y})}_{h}\right )-V(x,\phi)= \nonumber \\ && \qquad
\frac{\partial V(x,\phi)}{\partial x}z
+ \left .\frac {\partial V\left (x, \psi^{(x,\phi,{\cal
Y})}_{\ell}\right )}{\partial \ell}\right |_{\ell=0}h+o\left
(\sqrt{\tau^2+\vert z\vert^2+h^2}\right),\nonumber \\ &&
\end{eqnarray}
with $lim_{s\to 0^+} \frac {o\left ( \sqrt{s}\right )}{ \sqrt{s}
}=0$.
\end{definition}

 In the following, ${\mathcal V}$
is the class of functionals $V:R^n\times{\mathcal Q}\to R^+$ which
have the following properties: i) $V$ is locally Lipschitz in
${\mathcal C}$ and invariantly differentiable; ii) the maps, for
$\phi\in {\mathcal C}$ (involved $x\in R^n$, $h\in [0,\Delta)$),
\begin{eqnarray}&&
\phi\to \left .\frac {\partial V\left ( \phi(0),
\phi^h_{[-\Delta,0)}\right )}{\partial h}\right |_{h=0},
\qquad \phi \to \left . \frac{\partial
V(x,\phi_{[-\Delta,0)})}{\partial x}\right |_{x=\phi(0)}
\end{eqnarray}
are completely continuous and locally Lipschitz in ${\mathcal C}$.

\section {Main Results}
Let us consider a system described by the following RFDE
\begin{eqnarray}\label{equazionedisturbata}&&
\dot x(t)=f(x_t)+g(x_t)(u(t)+d(t)), \qquad t\ge 0,  \nonumber
\\ && x(\tau)=\xi_0(\tau),\qquad \tau \in [-\Delta,0], \qquad   \xi_0\in {\mathcal C},
\end{eqnarray}
which is exactly the same as (\ref{equazionediscontinua}), with
the same hypotheses, with the input signal given by the control
signal $u(t)\in R^m$ plus the actuator disturbance $d(t)\in R^m$,
Lebesgue measurable and locally essentially bounded.

For a given functional $V: R^n\times{\mathcal Q}\to R^+$ in the
class ${\mathcal V}$, let the maps $a: {\mathcal C}\to R^+$,
$b:{\mathcal C}\to R^{m}$ (row vector), $k:{\mathcal C}\to R^m$ be
defined as follows, for $\phi\in {\mathcal C}$ (involved $x\in
R^n$, $l\in [0,\Delta)$),

\begin{eqnarray}&& \label{controllaw} a(\phi)= \left . \frac{\partial V(x,\phi_{[-\Delta,0)})}{\partial x}\right |_{x=\phi(0)}f(\phi)
+ \left .\frac {\partial V\left (\phi(0),
\phi^{\ell}_{[-\Delta,0)}\right )}{\partial \ell}\right
|_{\ell=0}, \nonumber \\ && b(\phi)= \left .\frac{\partial
V(x,\phi_{[-\Delta,0)})}{\partial x}\right
|_{x=\phi(0)}g(\phi), \nonumber \\ && k(\phi)=
\left \{
\begin{array}{cc}-\frac{a(\phi)+\sqrt{a^2(\phi)+\vert
b(\phi)\vert^4}}{\vert b(\phi)\vert^2}b^T(\phi), & b(\phi)\ne 0
\\ \\ 0, & b(\phi)=0
\end{array}
\right .
\end{eqnarray}
For a positive real $r$, let $k_r:{\mathcal C}\to R^m$ be defined
as follows, for $\phi\in {\mathcal C}$,
\begin{equation}
k_r(\phi)= \left \{
\begin{array}{c}-\frac{a(\phi)+\sqrt{a^2(\phi)+\vert
b(\phi)\vert^4}}{\vert b(\phi)\vert^2}b^T(\phi), 
\qquad \qquad  \vert b(\phi)\vert > r, \\ \\
\\ -\frac{a(\phi)+\sqrt{a^2(\phi)+\vert
b(\phi)\vert^4}}{r^2}b^T(\phi), 
\qquad \qquad \vert
b(\phi)\vert \le r
\end{array}
\right .
\end{equation}

The following hypothesis will be used in the forthcoming theorem.
\begin{hypothesis}\label{standassoas} There exist a
functional $V: R^n\times {\mathcal Q}\to R^+$ in the class
${\mathcal V}$, with corresponding maps $a$, $b$, functions
$\alpha_1$, $\alpha_2$, $\alpha_3$ of class ${\mathcal
K}_{\infty}$, positive reals $r$, $p$ such that, $\forall \
\phi\in {\mathcal C}$:

\begin{itemize}
\item[i)]$\alpha_1(\vert \phi(0)\vert)\le V(\phi(0),\phi_{[-\Delta,0)})\le
\alpha_2(M_a(\phi))$;

\medskip

\item[ii)] $b(\phi)=0\ =>\ a(\phi)\le 0$; 

\medskip

   \item[iii)]$ a^2(\phi)+\vert b(\phi)\vert^4\ge \alpha^2_3(M_a(\phi)
    )$;

\medskip

\item[iv)] $\sup_{\{\psi\in {\mathcal C},\ 0<\vert b(\psi)\vert\le r\}}\ \frac{a(\psi)}{\vert
b(\psi)\vert}\le p$.

\end{itemize}
\bigskip
\end{hypothesis}

\bigskip

\begin{remark} The point (i) in Hypothesis \ref{standassoas} is standard in the ISS theory for systems described by RFDEs (see, for instance, \cite{PepeJiangSCL}). The point (ii)
is the standard key point in the theory of Control Liapunov
Functions (see, for instance, \cite{SontagUniversal}) and Control
Liapunov Functionals (see \cite{KarJiangESAIM}). The point (iii)
allows that the stabilizer obtained with the Sontag's formula is
such that the derivative in Driver's form of the functional $V$
satisfies a standard inequality for ISS concerns (see, for instance,
\cite{PepeJiangSCL}). The point (iv) is a key issue
in order to avoid problems related to non Lipschitz map $k$ (i.e.,
the Sontag's stabilizer). A similar condition was introduced in
\cite{Jankovic} in the framework of control Liapunov-Razumikhin
functions (see Assumption 1 in \cite{Jankovic}).
\end{remark}

\begin{theorem}\label{thiss} Let Hypothesis \ref{standassoas} be
satisfied. Then:

\begin{itemize}
\item[1)] the map $k_r:{\mathcal C}\to R^m$ is completely continuous and locally Lipschitz in ${\mathcal C}$;

\item[2)] there exists a function $\beta$ of class ${\mathcal
{KL}}$ and a function $\gamma$ of class ${\mathcal K}$, both
independent of $r$ and $p$, such that, chosen any positive real
$q$, for the closed loop system (\ref{equazionedisturbata}) with

\begin{equation}\label{univstabilizer}
u(t)=k_r(x_t)-qb^T(x_t),
\end{equation}
the solution exists for all $t\ge 0$ and, furthermore, satisfies
the following inequality

\begin{equation}\label{ineqquasiiss}
\vert x(t)\vert \le \beta(\Vert x_0\Vert_{\infty},t)+\gamma\left (
\sqrt{\frac{2}{q}}\Vert d_{[0,t)}\Vert_{\infty} \right
)+\gamma\left( \sqrt{\frac{2}{q}}(2p+r)\right )
\end{equation}
\end{itemize}
\end{theorem}

\begin{remark}
From the main Theorem in \cite{SontagPrincipe} (see also Theorem
3.1 in \cite{PepeJiangSCL}) it follows that the function $\gamma$
in Theorem \ref{thiss} is given, for $s\ge 0$, by
\begin{equation}
\gamma(s)=\alpha_1^{-1}\circ \alpha_2\circ \alpha_3^{-1}(s^2)
\end{equation}
\end{remark}

\begin{remark}\label{casogas}
Because of the inequality (\ref{ineqquasiiss}), the closed-loop
system (\ref{equazionedisturbata}), (\ref{univstabilizer}) is ISpS
(see Definition 2.1 in \cite{JiangMareelsWang}) with respect to
the disturbance $d(t)$. Notice in (\ref{ineqquasiiss}) that, if
the disturbance is bounded, the solution can achieve an
arbitrarily fixed small neighborhood of the origin by increasing
the control tuning parameter $q$.
\end{remark}
\begin{remark}
In the case the disturbance $d(t)\equiv 0$ and the map $k$ is
completely continuous and locally Lipschitz in ${\mathcal C}$,
then the control law

\begin{equation}\label{contgas}
u(t)=k(x_t)
\end{equation}
yields that the closed-loop system (\ref{equazionedisturbata})
(with $d(t)\equiv 0$), (\ref{contgas}) is $0-$GAS. In this case,
by Theorem 2.1, p. 132, in \cite{HaleLunel}, it is sufficient that
the point (i) of Hypothesis \ref{standassoas} holds with $\Vert
\phi\Vert_{\infty}$ instead of $M_a(\phi)$, $\phi \in {\mathcal
C}$, and the point (iii) holds with $\vert \phi(0)\vert$ instead
of $M_a(\phi)$, $\phi \in {\mathcal C}$. However, the locally
Lipschitz property of the map $k$ is a hypothesis which in general
is either not satisfied, either difficult to check. The
methodology here proposed treats also the actuator disturbance and
avoids the problem of the local Lipschitz property of the map $k$,
i.e. of the Arstein-Sontag stabilizer here extended to retarded
systems by the use of invariantly differentiable functionals. This
is done by the use of the map $k_r$, which is a slight
modification of the map $k$, and by the Sontag's ISS redesign
method (i.e., by adding the term $-qb$ in the control law
(\ref{controllaw})). The price to pay is that, in the case of
actuator disturbance $d(t)\equiv 0$, which in general is an
unrealistic hypothesis, an arbitrarily small neighborhood of the
origin is guaranteed to be asymptotically reached, instead of the
origin itself.
\end{remark}

\begin{remark}
The term $p$ in the point (iv) of Hypothesis \ref{standassoas} in
general can be reduced to $0$ by reducing $r$. Consequently, the
term $2p+r$ in (\ref{ineqquasiiss}) can be reduced to $0$. 
\end{remark}

\begin{remark} Notice that, as reported in the proof of Theorem \ref{thiss},
since the maps $k_r$ and $b$ are completely continuous and locally
Lipschitz in ${\mathcal C}$, it readily follows that the feedback
control map $k_r-qb^T$ is completely continuous and locally
Lipschitz in ${\mathcal C}$ as well.
\end{remark}

\begin{remark}
Notice that the methodology here presented, based on the
Arstein-Sontag approach with invariantly differentiable
Lyapunov-Krasovskii functionals, can be applied to time invariant,
control affine, nonlinear systems with an arbitrary number of
discrete as well as distributed time-delays, of arbitrary size.
The limitations introduced in \cite{Jankovic} to the form of the
maps $f$, $g$, describing the dynamics of the system (see (3.2) in
\cite{Jankovic}), are not necessary here.
\end{remark}

\section{Conlcusions}

In this paper we have considered an input-to-state practical
stabilizer, with respect to disturbances adding to the control
law, for retarded, control affine, nonlinear systems. We have used
the Sontag's formula, generated by means of control
Liapunov-Krasovskii functionals which are invariantly
differentiable. The results are very general and avoid problems
due to non Lipschitz control law at suitable subsets of the
infinite dimensional state space. This goal, together with
disturbance ISpS attenuation, is obtained by combining the
Sontag's formula, revisited in the above critical subsets, and the
Sontag ISS feedback control redesign method. If the actuator
disturbance is bounded, then an arbitrarily small neighborhood of
the origin
can be reached 
by increasing a
suitable control tuning parameter.


\begin{thebibliography}{2}
%
\def\IEEETAC{{\itshape IEEE Trans.~Automat.\ Contr.}}
\def\IEEETCASI{{\it IEEE Trans.~Circuits \&\ Sys., I}}
\def\IJC{{\itshape Int.~J.\ Contr.\/}}
\def\SCL{{\itshape Syst.\ Contr.\ Lett.\/}}
\def\SIAMJCO{{\itshape SIAM J.~Contr.\ Optim.\/}}


\bibitem{kybGMP} A. Germani, C. Manes, and P. Pepe, ``Local asymptotic stability for nonlinear state feedback delay systems,'' Kybernetika, Vol. 36, No. 1, pp. 31-42,
2000.

\bibitem{GerManPeIJRNC} A. Germani, C. Manes, and P. Pepe, ``Input-output
linearization with delay cancellation for nonlinear delay systems:
the problem of the internal stability,'' {\it International
Journal of Robust and Nonlinear Control}, Vol. 13, No. 9, pp.
909--937, 2003.



 \bibitem{HaleLunel}J.~K. Hale, and S.~M.~Verduyn Lunel, {\it
 Introduction to Functional Differential Equations}, Springer
 Verlag, 1993.



\bibitem{HuaGuanShi} C. Hua, X. Guan, and P. Shi, ``Robust stabilization of a
class of nonlinear time-delay systems,'' {\it Applied Mathematics
and Computation}, Vol.~155,
 pp.~737--752, 2004.



\bibitem{JiangMareelsWang} Z.-P. Jiang, I.M.Y. Mareels, Y. Wang, ``A Lyapunov formulation of the nonlinear small-gain theorem for
interconnected ISS systems,'' {\it Automatica}, Vol. 32, Issue 8,
pp. 1211-1215, 1996.

\bibitem{Jankovic} M. Jankovic, ``Control Lyapunov-Razumikhin functions
and robust stabilization of time delay systems,'' {\it IEEE
Transactions on Automatic Control}, Vol. 46, No. 7, pp.
1048--1060, 2001.

\bibitem{JankovicSepulchreKokotovic} M. Jankovic, R. Sepulchre, P.V. Kokotovic ``CLF based designs with robustness to dynamic input uncertainties,'' {\it Systems \& Control Letters},
Vol. 37, pp. 45--54, 1999.


\bibitem{KarJiangESAIM} I. Karafyllis, Z.-P. Jiang, ``Necessary and
Sufficient Lyapunov-Like Conditions for Robust Nonlinear
Stabilization'', {\it ESAIM: Control, Optimisation and Calculus of
Variations}, Vol. 16, pp. 887-928, 2010.

\bibitem{KarafyllisJiangBook} I. Karafyllis, Z.P. Jiang, {\it Stability and
Stabilization of Nonlinear Systems}, Springer, London, 2011.



\bibitem{KhalilBook} H. K. Khalil, {\it Nonlinear Systems}, Prentice Hall, International Edition, Third Edition, Upper Saddle River, New Jersey, 2000.

\bibitem{Kim}A.V. Kim, ``On the Lyapunov's functionals method for systems
with delays,'' {\it Nonlinear Analysis, Theory, Methods \&
Applications}, Vol. 28, No. 4, pp. 673-687, 1997.

\bibitem{KimBook}A.V. Kim, {\it Functional differential equations,
application of i-smooth calculus}, Kluwer Academic Publishers,
Dordrecht, 1999.


\bibitem{NEWKolmanovskii}V. Kolmanovskii and A. Myshkis, {\it Introduction
to the theory and applications of functional differential
equations}, Kluwer Academic Publishers, Dordrecht, 1999.









\bibitem{Lien}C.-H. Lien, ``Global exponential stabilization for several
classes of uncertain nonlinear systems with time-varying delay,''
{\it Nonlinear Dynamics and Systems Theory}, Vol.~4, No.~1,
 pp.~15--30, 2004.



\bibitem{MahboobiEsfanjaniNikravesh} R. Mahboobi Esfanjani, S.K.Y.
Nikravesh, ``Stabilising predictive control of non-linear
time-delay systems using control Lyapunov-Krasovskii
functionals'', {\it IET Control Theory and Applications}, Vol. 3,
Iss. 10, pp. 1395-1400, 2009.




\bibitem{MarquezMoog} L.A. Marquez-Martinez, and C.H. Moog, ``Input-output
feedback linearization of time-delay systems,'' {\it IEEE
Transactions on Automatic Control}, Vol.~49, No.~5, pp.~781-785,
2004.



\bibitem{OguchiIJC}  T. Oguchi, A. Watanabe and T. Nakamizo,
``Input-output linearization of retarded non-linear systems by
using an extension of Lie derivative,'' {\it International Journal
of Control}, Vol.~75, No.~8, 582--590, 2002.




\bibitem{PepeIJACSP} P. Pepe, ``Adaptive Output Tracking for a Class of Nonlinear Time Delay Systems,'' {\it International Journal of Adaptive Control and Signal Processing}, Vol. 18, No. 6,
pp. 489-503, 2004.




\bibitem{PepeTAC09}
P. Pepe, ``Input-to-state stabilization of stabilizable,
time-delay, control affine, nonlinear systems,'' {\it IEEE
Transactions on Automatic Control}, Vol. 54, No. 7, pp. 1688-1693,
2009.


\bibitem{PepeJiangSCL}P. Pepe, and Z.-P. Jiang, ``A Lyapunov-Krasovskii
methodology for ISS and iISS of time-delay systems,'' {\it Systems
\& Control Letters}, Vol.~55, No.~12,
 pp.~1006--1014, 2006.



\bibitem{SontagPrincipe}E.D. Sontag, ``Smooth stabilization implies
 coprime factorization, {\it IEEE Transactions on Automatic
 Control,} Vol. 34, No. 4, pp. 435--443, 1989.

 \bibitem{SontagUniversal}E.D. Sontag. ``A "universal" construction of Artstein's theorem on nonlinear stabilization'', {\it Systems \& Control Letters}, Vol. 13, No. 2, pp. 117-123, 1989.


 \bibitem{ZhangCheng}X. Zhang, Z. Cheng, ``Global stabilization of a class
 of time-delay nonlinear systems,'' {\it International Journal of
 Systems Science}, Vol.~36, No.~8,
 pp.~461--468, 2005.



\end{thebibliography}
\end{document}